\documentclass{amsart}
\usepackage{amssymb, epic,eepic,epsfig,amsbsy,amsmath,amscd}

\def\be#1\ee{\begin{equation}#1\end{equation}}
\theoremstyle{plain}

\newtheorem{theorem}{Theorem}
\newtheorem{proposition}{Proposition}[section]
\newtheorem{lemma}[proposition]{Lemma}
\newtheorem{corollary}[proposition]{Corollary}

\theoremstyle{definition}

\theoremstyle{remark}

\def\printname#1{
    \if\draft y
        \smash{\makebox[0pt]{\hspace{-0.5in}
            \raisebox{8pt}{\tt\tiny #1}}}
    \fi
}

\newlength{\standardunitlength}
\newcommand{\qbinom}[2]{\text{$\left[\begin{array}{c}#1\\ #2\end{array}
\right]$}}

\def\BZ{\mathbb Z}

\def\BQ{\mathbb Q}

\def\BC{\mathbb C}

\def\cL{\mathcal L}

\def\N{{\mathbb N}}

\def\R{\mathcal R}

\newcommand{\sn}{\operatorname{sn}}

\def\Z{\BZ}

\def\msl{\mathfrak{sl}}

\def\ev{\mathrm{ev}}

\newcommand{\Gd} {\Lambda_a}
\newcommand{\hGd} {\hat \Lambda_a}

\def\Habiro{\widehat{\Z[q]}}

\newcommand{\psdiag}[3]{\hspace{1mm}\raisebox{-#1mm}{\epsfysize#2mm
\epsffile{#3.eps}}\hspace{1mm}}

\begin{document}

\title[Integrality $\&$ Rational Surgery]{ Integrality
of quantum 3--manifold invariants\\ and rational
surgery formula}

\author{Anna Beliakova}
\address{Institut f\"ur Mathematik, Universit\"at Zurich,
Winterthurerstrasse 190, 8057 Z\"urich, Switzerland}
\email{anna@math.unizh.ch}

\author[Thang Le]{Thang T. Q. L\^e}
\address{Department of Mathematics \\
         Georgia Institute of Technology \\
         Atlanta, GA 30332-0160, USA}
\email{letu@math.gatech.edu}

\begin{abstract}
We prove that the Witten--Reshetikhin--Turaev (WRT) $SO(3)$ invariant
of an arbitrary 3--manifold $M$ is always an algebraic integer.
Moreover, we give a rational surgery formula for
 the unified invariant
 dominating  WRT $SO(3)$ invariants of rational homology
3--spheres at roots of unity of order co--prime with the torsion.
 As an application,
we compute the unified
invariant for Seifert fibered spaces and for Dehn surgeries on twist knots.
We show that this invariant separates
integral homology Seifert fibered spaces and can be used to detect
the unknot.
\end{abstract}

\maketitle

\section*{Introduction}

The Witten--Reshetikhin--Turaev (WRT) invariant was first introduced
by Witten using physics heuristic ideas, and then mathematically
rigorously by Reshetikhin and Turaev \cite{Tur}. The invariant, depending on a
root $\xi$ of unity, was first defined for the Lie group $SU(2)$, and
was later generalized to other Lie groups. The $SO(3)$ version of
the invariant was introduced by Kirby and Melvin \cite{KM1}. For
this
 $SO(3)$ version the quantum parameter $\xi$ 
must be a root of unity of {\em odd} order.
 One important result in quantum topology, first proved by H.
Murakami for rational homology spheres \cite{Murakami} and then
 generalized  by Masbaum and Roberts \cite{MR},
is that the WRT $SO(3)$ invariant (also known as quantum
$SO(3)$ invariant)
 $\tau_M(\xi)$ of an arbitrary 3--manifold $M$
is an algebraic integer, when the order of the root of unity $\xi$
is an {\em odd prime}. Recently,  the second author proved
\cite{Le} that if the order of  $\xi$  is co--prime with the
cardinality of the torsion of $H_1(M,\BZ)$, then
 the $SO(3)$ quantum invariant $\tau_M(\xi)\in \BZ[\xi]$.
 In this paper we remove all the restrictions  on the order of $\xi$.
\begin{theorem} \label{main} For every closed 3--manifold and
 every root $\xi$ of unity of odd order, the
quantum $SO(3)$ invariant $\tau_M(\xi)$  belongs to $\BZ[\xi]$.
\end{theorem}

The integrality has many important applications, among them is the
construction  of an integral topological quantum field theory
 and representations of mapping class groups
over $\BZ$ by Gilmer and  Masbaum (see e.g. \cite{GM}).
The integrality is also a
 key property required for the categorification
of  quantum 3--manifold invariants \cite{K}.

Our proof of integrality is inspired by the Habiro's work.
In  \cite{Hab}, Habiro constructed
an invariant of integral homology 3--spheres
with values in the universal ring
$$\Habiro=\lim_{\overleftarrow{\hspace{2mm}n\hspace{2mm}}}
\frac{ \Z[q]}{ (q)_n}\, $$ where
$(q)_n:=(q;q)_n=(1-q)(1-q^2)\dots(1-q^n)$. Habiro's invariant
specializes at a root $\xi$ of unity to $\tau_M(\xi)$.

In \cite{Le}, the second author generalized  Habiro's theory to
rational homology 3--spheres. For a rational homology sphere $M$
with $|H_1(M,\Z)|=a$, he constructed an invariant $I_M$ which
dominates the $SO(3)$ invariants of $M$ at roots of unity of order
co--prime to $a$. Habiro's universal ring was modified by inverting
$a$ and cyclotomic polynomials of order not co--prime to $a$.
Applications of this theory are the new integrality properties of
quantum invariants, new results about Ohtsuki series and a  better
understanding of the relation between LMO invariant, Ohsuki series
and  quantum invariants.

In  this paper we give a rational surgery formula for
the unified invariant $I_M$ defined in \cite{Le}.
The main applications of our construction are summarized below.

For a positive integer $a$  let $ A_a := \BZ[\frac{1}{a}][q^{\pm
1/a}]$ and $\N_a$ the set of positive integers co--prime
with $a$.
 Denote by $\Phi_s(t)$ the $s$--th cyclotomic polynomial. Let
$\Gd \subset \BQ(q^{1/a})$ be the ring obtained from $A_a$ by adding
the inverses of each $\Phi_s(q^{1/a})$ with $s$ {\em not co--prime}
with $a$:
$$\Gd := A_a\left[\frac{1}{\Phi_{s}(q^{1/a})}, \, s\notin \N_a\, \right]$$
The analog of the Habiro ring constructed in \cite{Le} is
$$\hGd := \lim_{\overleftarrow{\hspace{2mm}n\hspace{2mm}}}
\frac{\Gd}{(q)_n}$$
Let $U_a$  be the set of all complex roots of unity with orders {\em
odd} and {co--prime} with $a$. It was shown
in \cite{Le}, that $I_M\in \hGd$ dominates $SO(3)$ quantum invariants
of $M$ with $|H_1(M,\BZ)|=a$ at roots of unity from $U_a$.

Let $t:=q^{1/a}$. Let $\R_{a,k} \subset \Gd$ be a subring generated
over $\Z[t^{\pm 1}]$ by
  $\frac{(t;t)_{k}}{(q;q)_{k}}$. Note that,

  $$ \R_{a,1} \subset \R_{a,2} \subset \dots \subset \R_a,$$

Where $\R_a= \cup_{k=1}^\infty \R_{a,k}$, also a subring of
$\Lambda_a$. Unlike the case of $\Lambda_a$, in the construction of
  $\R_a$ we don't need to invert $a$.
Let
$$\widehat \R_a := \lim_{\overleftarrow{\hspace{2mm}n\hspace{2mm}}}
\frac{\R_a}{(q)_n}$$ be its cyclotomic completion. Then we can
refine the result of \cite{Le} as follows.
\begin{theorem}\label{ring} Let $M$ be a rational homology 3--sphere.
Then there exist $m\in \N$ and
 $f_{k_i}(t)\in \R_{a,2k_{i}+1}$ for $i=1,\dots,m$, such that
 the unified  invariant $I_M$ admits the following presentation.
$$I_M=\sum^\infty_{k=0}\; \sum^k_{k_1,...,k_m\geq 0}\;
\prod^m_{i=1} \;f_{k_i}(t) \;(q^{k+1})_{k+1}.$$

In particular, $I_M \in \widehat \R_a$.

\end{theorem}

Further, we compute the unified invariant for Seifert fibered spaces
and for Dehn surgeries on twist knots.
\begin{theorem}\label{SFS}
The unified invariant  separates integral homology Seifert fibered spaces.
\end{theorem}

For a knot $K$, let $M(K,a)$ denotes the 3--manifold obtained by
surgery on the knot $K$ with framing $a$. In general, there are
different $K,K'$ such that $M(K,a)=M(K',a)$ for some $a$.

\begin{theorem}\label{detector}
Suppose $I_{M(K,a)}=I_{M(K',a)}$ for infinitely many $a \in \Z$.
Then $K$ and $K'$ have the same colored Jones polynomial.
\end{theorem}

In particular, using the recent deep result of Andersen \cite{JA},
that the colored Jones polynomial detects the unknot, we see that
(under the assumption of the theorem), if $K$ is the unknot, then so
is $K'$.

\subsection*{Acknowledgment}  The first author wishes to express her
gratitude to  Christian Krattenthaler for
valuable suggestions.

\section{Quantum invariants}
\label{defs}
Let us first fix the notation.

$$ \{n\} = q^{n/2}-q^{-n/2},
 \quad  \{n\}!=
\prod_{i=1}^n \{i\} ,\quad  [n] =\frac{\{n\}}{\{1\}}, \quad
\qbinom{n}{k} = \frac{\{n\}!}{\{k\}!\{n-k\}!}.$$

\subsection{The colored Jones polynomial}

\newcommand{\RR} {\mathbf R}

 Suppose $L$ is framed, oriented link
in $S^3$ with $m$ ordered components.
 For every positive integer $n$ there is a unique
irreducible $sl_2$--module $V_n$ of dimension $n$.
For positive integers $n_1,\dots,n_m$ one can define
the quantum invariant $J_L(n_1,\dots,n_m):=
J_L(V_{n_1},\dots,V_{n_m})$ known as the colored
Jones polynomial of $L$ (see e.g. \cite{Tur}).
Let us recall here a few well--known formulas.
For the unknot $U$ with 0 framing one has
\begin{equation} J_U(n) = [n]= \{n\}/\{1\}. \label{unknot}
\end{equation}
 If $L'$ is obtained from $L$
by increasing the framing of the $i$--th component by 1, then
\begin{equation}\label{framing}
J_{L'}(n_1,\dots,n_m) = q^{(n_i^2-1)/4} J_{L}(n_1,\dots,n_m).
\end{equation}
In general, $J_{L}(n_1,\dots,n_m) \in \BZ[q^{\pm 1/4}]$. However,
there is a number $a\in \{0,\frac{1}{4},\frac{1}{2},\frac{3}{4}\}$
such that $J_{L}(n_1,\dots,n_m) \in q^a\BZ[q^{\pm 1}]$.

\subsection{Evaluation map and Gauss sum}
Throughout this paper let $\xi$ be  a primitive root
of unity of {\em odd} order $r$.
We first define, for each $\xi$, the evaluation map
$\ev_\xi$, which replaces $q$ by $\xi$. Suppose $f\in \BQ[q^{\pm
1/h}]$, where $h$ is co--prime with $r$, the order of $\xi$. There
exists an integer $b$, unique modulo $r$, such that
$(\xi^b)^{h}=\xi$. Then we define
$$\ev_\xi f := f|_{q^{1/h}= \xi^b}.$$
The definition extends to $\ev_\xi: \hGd \to \BC$, since
$\ev_\xi((q;q)_n)=0$ if $n \ge r$.

Suppose
  $f(q;n_1,\dots,n_m)$ is a function
of variables $q$ and integers $n_1,\dots,n_m$. Let
$$ {\sum_{n_i}}^\xi f := \sum_{n_i} \ev_\xi (f),$$
where in the sum  all the $n_i$ run the set of {\em odd} numbers
between $0$ and $2r$.
 A variation
$\gamma_d(\xi)$ of the Gauss sum is defined by
$$ \gamma_d(\xi):= {\sum_{n}}^\xi q^{d\frac{n^2-1}{4}}.$$
It is known that, for odd $r$,
 $|\gamma_d(\xi)|= \sqrt r$, and hence  is never 0.
$$ F_L(\xi):= {\sum_{n_i}}^\xi
J_L(n_1,\dots,n_m)\prod_{i=1}^m [n_i].$$
 The following result is well--known
(compare  \cite{Le}).
\begin{lemma} For the unknot $U^\pm$ with framing $\pm 1$, one has
$F_{U^{\pm}}(\xi) \neq 0$. Moreover,
\be\label{unpm}
F_{U^{\pm}}(\xi) = \mp 2  \gamma_{\pm 1}(\xi)\, \ev_\xi\left( \frac{
q^{\mp 1/2}}{\{1\}}\right).
\ee
\label{un1}
\end{lemma}

\subsection{Definition of $SO(3)$ invariant of 3--manifolds}\label{def}
All 3--manifolds in this paper are supposed to be closed and
oriented. Every link in a 3--manifold is framed, oriented, and has
components ordered.

 Suppose $M$ is an
oriented 3--manifold obtained from $S^3$ by surgery along a framed,
oriented link $L$. (Note that $M$ does not depend on the orientation
of $L$). Let $\sigma_+ $ (respectively, $\sigma_-$) be the number of
positive (resp. negative) eigenvalues of the linking matrix of $L$.
Suppose $\xi$ is a root of unity of odd order $r$. Then the quantum
$SO(3)$ invariant is defined by
\begin{equation*}
\tau_M(\xi) = \tau_M^{SO(3)}(\xi) :=
\frac{F_L(\xi)}{(F_{U^{+}}(\xi))^{\sigma_+}\,
(F_{U^{-}}(\xi))^{\sigma_-} }.
\end{equation*}
For connected sum, one has $ \tau_{M\#N}(\xi) =\tau_{M}(\xi)
\tau_N(\xi).$

\subsection{Laplace transform}
In \cite{BBL}, we together with Blanchet developed the Laplace transform method
to compute $\tau_M(\xi)$.  Here we generalize this method to the case
where $r$ is not co--prime with torsion.

Suppose $r$ is an odd number, and $d$ is positive integer. Let
$$ c:= (r,d), \quad d_1:= d/c, \quad r_1:=r/c.$$
\newcommand{\Tor}{{\rm Tor}}

Let $\cL_{d;n}: \BZ[q^{\pm n},q^{\pm 1}] \to \BZ[q^{\pm 1/d}]$ be
the $\BZ[q^{\pm 1}]$-linear operator, called the Laplace transform,
defined by
\be\label{four}
\cL_{d;n}(q^{na}) := \begin{cases}
0
&\text{if $c\not| a$;}\\
q^{-a^2/d}   & \text{if $a=ca_1$},
 \end{cases}
\ee
\begin{lemma} Suppose  $f \in \BZ[q^{\pm n},q^{\pm 1}]$. Then
$$ {\sum_{n}}^\xi q^{d\frac{n^2-1}{4}}  f = \gamma_d(\xi) \,
\ev_\xi(\cL_{d;n}(f)).
$$
\label{1000}
\end{lemma}

\begin{proof} It's enough to consider the case when $f= q^{na}$,
with $a$ an integer. This case is proven by Lemma \ref{33}
in the next subsection.
\end{proof}

The point is that $\cL_{d;n}(f)$, unlike the left hand side
${\sum_{n}^\xi}  q^{d\frac{n^2-1}{4}} f$, does not depend on $\xi$,
and will help us to define a ``universal invariant".
Note that
Lemma \ref{1000} with $d=\pm1$ and $f= [n]^2$ implies Lemma \ref{un1}.

\subsection{Reduction from $r$ to $r_1$} Let $O_r$ be the set of all
odd integers between $0$ and $2r$. This set $O_r$ can be partitioned
into $r_1$ subsets $O_{r;s}$ with $s \in O_{r_1}$,  where $O_{r;s}$
is the set of all $n \in O_r$ which are equal to $s$ modulo $r_1$. In
other words,  $O_{r;s}:= \{ s + 2jr_1, j=0,1, \dots, c-1\}$. The
point is, the value of $\xi^{d \frac{n^2-1}{4}}$ remains the same
for all $n$ in the same set $O_{r;s}$.
Let $\zeta = \xi^c$, then the order of $\zeta$ is $r_1$.
\begin{lemma}
\label{33}
One has
\begin{equation}
 \gamma_d(\xi) = c \gamma_{d_1}(\zeta). \label{02}
 \end{equation}
\begin{equation}{\sum_{n}}^\xi q^{d \frac{n^2-1}{4}}\, q^{an} = \begin{cases}
0
&\text{if $c\not | a$;}\\
 \gamma_d(\xi) \zeta^{-a_1^2 d_1^*}  & \text{if $a=ca_1$},
 \end{cases}
 \end{equation}
where $d_1$ is an integer satisfying $d_1 d_1^* \equiv 1 \pmod{r_1}$.
\end{lemma}

\begin{proof}  One has
 \begin{align*} {\sum_{n}}^\xi q^{d \frac{n^2-1}{4}}\, q^{an} = \sum_{n\in O_r} \, \xi^{d \frac{n^2-1}{4}}\,
 \xi^{an}
 = \sum_{s\in O_{r_1}}  \sum_{n\in O_{r;s}} \, \xi^{d
 \frac{n^2-1}{4}} \xi^{an}.
 \end{align*}
Using the fact that   $\xi^{d \frac{n^2-1}{4}}$ remains the same for
all $n$ in the same set $O_{r;s}$, we get
\begin{align}
{\sum_{n}}^\xi q^{d \frac{n^2-1}{4}}\, q^{an} & = \sum_{s\in
O_{r_1}} \xi^{d
 \frac{s^2-1}{4}} \sum_{n\in O_{r;s}} \xi^{an} \label{05}\\
 & = \sum_{s\in
O_{r_1}} \xi^{d \frac{s^2-1}{4}} \xi^{sa} \left( \sum_{j=0}^{c-1}
\xi^{2ar_1j} \right) \label{06}
\end{align}
Note that \eqref{02} follows from \eqref{05} with $a=0$.
\begin{equation}
 \sum_{j=0}^{c-1} \xi^{2ar_1j} = \sum_{j=0}^{c-1}
(\xi^{2ar_1})^j. \label{01}
\end{equation}
If $c\not| a$, then $(\xi^{2ar_1}) \neq 1$, but a root of unity of
order dividing $c$, hence the right hand side of \eqref{01} is 0.
It follows that the right hand side of \eqref{06} is also 0, or
${\sum_{n}}^\xi q^{d \frac{n^2-1}{4}}\, q^{an} =0$.

If $c|a$, then the right hand side of \eqref{01} is $c$. Hence from
\eqref{06} we have

\begin{align*}
 {\sum_{n}}^\xi q^{d \frac{n^2-1}{4}}\, q^{an} & = c \sum_{s \in O_{r_1}}
\xi^{d \frac{s^2-1}{4}} \xi^{sa} \\
&= c \sum_{s \in O_{r_1}} \zeta^{d_1 \frac{s^2-1}{4}} \zeta^{sa_1} =
c {\sum_{n}}^\zeta q^{d_1
\frac{n^2-1}{4}}\, q^{a_1n}\\
&= c \gamma_{d_1}(\zeta) \zeta^{-a_1^2 d_1^*}  \qquad
\end{align*}
The last equality follows by the standard square completion argument.
Using  \eqref{02} we get the result.
\end{proof}

\subsection{Habiro's cyclotomic expansion of the colored Jones
polynomial}

\label{300}
In \cite{Hab},  Habiro defined a new basis
$P'_k$, $k=0,1,2,\dots$,  for the Grothendieck ring of finite--dimensional
$sl_2$--modules, where
$$ P_k':=
\frac{1}{\{k\}!} \, \prod_{i=1}^{k}(V_2-q^{(2i-1)/2}-q^{-(2i-1)/2})
.
$$
For any  link $L$, one has
\begin{equation}\label{Jones}
J_{L}(n_1,\dots,n_m) = \sum_{0\le k_i \le n_i-1} J_L(P'_{k_1}, \dots
, P'_{k_m}) \prod_{i=1}^m \qbinom{n_i+k_i}{2k_i+1} \{k_i\}!
\end{equation}
Since there is a denominator in the definition of $P_k'$, one might
expect that $J_L(P'_{k_1}, \dots , P'_{k_m})$ also has non--trivial
denominator. A difficult and important integrality result of Habiro
\cite{Hab}  is

\begin{theorem}\cite[Thm. 3.3]{Hab}
\label{Hab} If $L$ is algebraically split and zero framed link in
$S^3$, then
$$
J_L(P'_{k_1}, \dots , P'_{k_m}) \in \frac{\{2k+1\}!}{\{k\}!\{1\}}
\,\,\BZ[q^{\pm1/2}] = \qbinom{2k+1}{k} (q^2)_k \,\BZ[q^{\pm1/2}],
$$
where $k=\max\{k_1,\dots, k_m\}$.
\end{theorem}
Thus, $J_L(P'_{k_1}, \dots , P'_{k_m})$ is not only integral, but
also divisible by $(q)_k$.

 Suppose
$L$ is an algebraically split link with 0--framing on each component.
Then we have
\begin{equation*}
\ev_\xi (J_{L}(n_1,\dots,n_m))  = \ev_\xi \left(
\sum_{k_1,\dots,k_m=0}^{(r-3)/2} J_L(P'_{k_1}, \dots , P'_{k_m})
\prod_{i=1}^m \qbinom{n_i+k_i}{2k_i+1} \{k_i\}!\right)
\end{equation*}

\section{Integrality of quantum invariants for all roots of unity}
Throughout this section we assume
that   $c=(d,r)>1$, $r/c=r_1$,
$d/c=d_1$ and $d_1d^*_1=1 \mod r_1$, where $r$ is the order of $\xi$
and $d$ is the order of the torsion part of $H_1(M,\BZ)$.

\subsection{Quantum invariants of links with diagonal linking matrix}
The following proposition
plays a key role in the proof of integrality.

\begin{proposition} \label{222}
 For  $k \le
 (r-3)/2$, we have
\be
\frac{1}{\gamma_{\pm 1}(\xi)} {\sum_{n}}^\xi  q^{d\frac{n^2-1}{4}} \,
\qbinom{n+k}{2k+1} \{k\}! \{n\} \; \in \;\BZ[\xi] . \label{121}
\ee
\end{proposition}

\vspace*{0.2cm}

\noindent
{\bf Proof of Theorem \ref{main} (diagonal case)}
Suppose $M$ is obtained from $S^3$ by surgery along an algebraically
split $m$--component link $L$ with integral framings $d_1, d_2, \dots,d_m$.
Inserting into the definition of $\tau_M(\xi)$ (see Section \ref{def})
the formulas \eqref{unpm}
and \eqref{Jones} and using Lemma \ref{1000}, we see
that Proposition \ref{222} and Theorem \ref{Hab} imply integrality
if  $d_i\neq 0$ for all $i$. If some of $d_i$ are zero, then
by  same argument as  in Section 3.4.2 of \cite{Le}  we have
$$
{\sum_{n}}^\xi   \qbinom{n+k}{2k+1} \{k\}! \{n\} ,\ = \, 2
\ev_\xi\left( q^{(k+1)(k+2)/4} \, (q^{k+2})_{r-k-2}
\right) \, .
$$
The result follows now from the fact
$\gamma_d(\xi)/\gamma_1(\xi)\in \BZ[\xi]$.
\qed

\subsubsection{Technical results} This subsection is devoted
to the proof of  Proposition \ref{222}.

\begin{lemma} {\rm (a)} Suppose $x \in \BQ(\xi)$ such that $x^2 \in
\BZ[\xi]$, then $x \in \BZ[\xi]$.

{\rm (b)} Suppose $x,y \in \BZ[\xi]$ such that $x^2$ is divisible by
$y^2$, then $x$ is divisible by $y$. \label{2001}
\end{lemma}

\begin{proof}(a) Suppose $a=x^2$, then $a\in \BZ[\xi]$ and $x$ is a
solution of $x^2-a=0$, hence $x$ is integral over $\BZ[\xi]$, which
is integrally closed. It follows that $x\in \BZ[\xi]$.

(b) We have that $(x/y)^2 = x^2/y^2$ is in $\BZ[\xi]$, hence by part
(a), $x/y\in \BZ[\xi]$.
\end{proof}

Recall that
$$(q^l;q)_m= \prod_{j=l}^{l+m-1}(1-q^j).$$

 Let
$\widetilde{(q^l;q)}_m$ be the product on the right hand side, only with
$j$ {\em not divisible by $c$}. Also let $\widehat{(q^l;q)}_m$ be the
complement, i.e. $\widehat{(q^l;q)}_m :=
(q^l;q)/\widetilde{(q^l;q)}_m.$
Using
$ (\xi;\xi)_{r-1} = r, \quad (\xi^c;\xi^c)_{r_1-1} = r_1,$
we see that
\be \widetilde{(\xi;\xi)}_{r-1} =c, \label{ccc}
 \ee
where $(a;b)_m:=(1-a)(1-ab)\dots(1-ab^{m-1})$.
Note that $1-\xi^j$ is invertible in $\BZ[\xi]$ iff $(j,r)=1$. Let
\be\label{z} z:=  \widetilde {(\xi;\xi)}_{(r-1)/2}, \quad \text
{and} \quad z' := \widetilde {(\xi^{(r+1)/2};\xi)}_{(r-1)/2}. \ee
Then $zz'$ is the left hand side of \eqref{ccc}, hence $zz'=c$.
We use the notation  $x\sim y$ if the ratio $x/y$ is a unit in
$\BZ[\xi]$. Note that $z\sim z'$. This is because
   $1-\xi^k \sim 1-\xi^{r-k}$. Thus we have

   \be z^2 \sim c\label{2002}
   \ee

\begin{lemma}\label{divisgamma}
 $\gamma_d(\xi)/\gamma_1(\xi)$ is divisible by $z$.
\end{lemma}
\begin{proof} Using Lemma \ref{2001}(b) and \eqref{2002}, one needs only to show that
$(\gamma_d(\xi))^2/(\gamma_1(\xi))^2$ is divisible by $c$.  The
values of $\gamma_b(\xi)$ are well--known when $b$ is co--prime with
$r$, the order of $\xi$. In particular, $\gamma_b(\xi) \sim
\gamma_1(\xi)$, see \cite{LL}.

Recall that $\zeta = \xi^c$ has order $r_1$. Since $d_1$ and $r_1$
are co--prime, we have
$$ \gamma_{d_1}(\zeta) \sim \gamma_{1}(\zeta).$$
Using Lemma \ref{33}, we have
\be \frac{(\gamma_{d}(\xi))^2}{(\gamma_{1}(\xi))^2} = c^2
\frac{(\gamma_{d_1}(\zeta))^2}{(\gamma_{1}(\xi))^2} \sim c^2
\frac{(\gamma_{1}(\zeta))^2}{(\gamma_{1}(\xi))^2} \label{2003}\ee
Using explicit formula for $\gamma_{1}(\xi)=\sum_{0\leq j<r}\xi^{j^2+j}$
(given e.g. by Thm. 2.2 of \cite{LL}),
we have that
$$ (\gamma_{1}(\xi))^2 = \pm r \xi^{-2^*} = cr_1 \xi^{-2^*} , \quad
(\gamma_{1}(\zeta))^2= \pm r_1 \zeta^{-2^*}$$
where $2^*$ is the inverse of $2$.
Plugging this in \eqref{2003}, we get the result.
\end{proof}

For $c, b\in \BZ$ we define
\be
Y_c(k,b)  :=  (-1)^k \sum_{n=-\lfloor k/c\rfloor}^{\lfloor(k+1)/c\rfloor}
 (-1)^n
\qbinom{2k+1}{k+nc}\,
q^{cbn^2}
\ee

\begin{lemma} \label{1110}
Suppose $d_1 d_1^* \equiv 1 \pmod {r_1}$,
where $r=cr_1$ is  the order of $\xi$,
then
\begin{align*} {\sum_{n}}^\xi  q^{d\frac{n^2-1}{4}} \,  \qbinom{n+k}{2k+1} \{k\}!
\{n\} & = -2 \gamma_d(\xi) \, \ev_\xi \left(\frac{\tilde
Y_c(k,-d^*_1)\{k\}!}{\{2k+1\}!}\right).
\end{align*}
\end{lemma}

\begin{proof}
By Lemma \ref{1000} we have to compute
 $\cL_{d;n}(\{n\} \, \{n+k\}!/\{n-k-1\}!)$.
 Since $\cL_{d;n}$ is invariant under $n \to -n$, one has
\begin{equation}\cL_{d;n}(\{n\} \, \{n+k\}!/\{n-k-1\}!) = -2 \cL_{d;n}(q^{-nk}\,(q^{n-k};q)_{2k+1}
).\label{123}
\end{equation}
By the $q$--binomial formula we have
\begin{equation} q^{-nk} (q^{n-k};q)_{2k+1} = \sum_{j=0}^{2k+1}(-1)^j
\qbinom{2k+1}{j} \, q^{n(j-k)}. \label{4}
\end{equation}
Using the definition of $\cL_{d;n}$  we get
\begin{align} \ev_\xi(
 \cL_{d;n}(\{n\} \, \{n+k\}!/\{n-k-1\}!))
&= -2 \ev_\xi (Y_c(k,-d^*_1)). \notag
\end{align}
Multiplying by $\{k\}!/\{2k+1\}!$, we get the result.
\end{proof}

\begin{theorem}\label{techx} For
$b\in \BZ$ and $k\leq (r-3)/2$,
 $\frac{\gamma_d(\xi)}{\gamma_1(\xi)}
\ev_{\xi}(Y_c(k,b))$ is divisible by $\ev_{\xi}\left(\frac{\{2k+1\}!}{\{k\}!}
\right)$.
\end{theorem}
Here we modify the proof of Theorem 7 in \cite{Le}.
\begin{proof}
The case $b=0$ is trivial. Let us assume $b\neq 0$.
Separating the case $n=0$ and combining
positive and negative $n$, we have
$$ Y_c(k,b) = (-1)^k  \qbinom{2k+1}{k} + (-1)^k
\sum_{n=1}^{\lfloor k/c\rfloor}(-1)^n q^{cb
n^2} \left(  \qbinom{2k+1}{k+nc} + \qbinom{2k+1}{k-nc} \right) .$$
Using
$$ \qbinom{2k+1}{k+cn} + \qbinom{2k+1}{k-cn} =
\frac{\{k+1\}}{\{2k+2\}}\qbinom{2k+2}{k+cn+1}(q^{cn/2} +
q^{-cn/2}) $$ and
$ \qbinom{2k+2}{k+1} =
\qbinom{2k+1}{k}\frac{\{2k+2\}}{\{k+1\}}$ we get

\begin{equation}
Y_c(k,b) = (-1)^{k} \qbinom{2k+1}{k} S_{N} \label{001x}
\end{equation}
where $N=k+1$ and
$$S_N= 1 + \sum_{n=1}^\infty \frac{q^{Ncn} (q^{-N};q)_{cn}}{(q^{N+1};q)_{cn}}\,
(1+q^{cn}) \, q^{cbn^2}.$$
For  $z$ defined
by \eqref{z}, we show
the divisibility of $\ev_{\xi} (S_N) z$ by $(\xi;\xi)_N$
in Section \ref{ident}. This
implies the result, since  $z |  \frac{\gamma_d(\xi)}{\gamma_1(\xi)}$
by Lemma \ref{divisgamma} and
$$\frac{\qbinom{2k+1}{k} \{k\}!}{\{2k+1\}!}=\{k+1\}!\, .$$

\end{proof}

\vspace*{0.2cm}
\noindent
{\bf Proof of Proposition \ref{222}}
Combining Lemma \ref{1110} with Theorem \ref{techx} we get
Proposition \ref{222}.
\qed
\subsubsection{Andrew's identity}\label{ident}

Let $\alpha_n,\beta_n$ be a Bailey pair as defined in Section 3.4 of
\cite{Andrews}, with $a=1$. Then for any numbers $b_i,c_i,
i=1,\dots,k$ and positive integer $N$ we have the identity (3.43) of
\cite{Andrews}:
\begin{multline}\sum_{n\ge 0} (-1)^n\alpha_n q^{-\binom{n}{2}+kn+Nn}
\frac{(q^{-N})_n}{(q^{N+1})_n} \prod_{i=1}^k
\frac{(b_i)_{n}}{b_i^{n}} \frac{(c_i)_{n}}{c_i^{n}}
\frac{1}{(\frac{q}{b_i})_n (\frac{q}{c_i})_n}=
\\
\frac{(q)_N\,  (\frac{q}{b_kc_k})_N}{(\frac{q}{b_k})_N\,
(\frac{q}{c_k})_N} \sum_{n_k \ge n_{k-1} \ge \dots \ge n_1\ge 0}
\beta_{n_1} \frac{q^{n_k} (q^{-N})_{n_k} (b_k)_{n_k}
(c_k)_{n_k}}{(q^{-N}b_kc_k)_{n_k}}\prod_{i=1}^{k-1}\frac{q^{n_i}
\frac{(b_i)_{n_i}}{b_i^{n_i}} \frac{(c_i)_{n_i}}{c_i^{n_i}}
(\frac{q}{b_i c_i})_{n_{i+1}-n_i} }{(q)_{n_{i+1}-n_i}
(\frac{q}{b_i})_{n_{i+1}}  (\frac{q}{c_i})_{n_{i+1}}}. \label{1002}
\end{multline}
A special Bailey pair is given by (see section 3.5 of
\cite{Andrews}):
\begin{align*} \alpha_0 &=1, \quad \alpha_n = (-1)^n q^{n(n-1)/2}(1+q^n) \quad
\text{for } n \ge 1.\\
\beta_0 & =1, \quad \beta_n =0 \quad \text{for } n \ge 1.
\end{align*}
 Using the decomposition
$$(q^x;q)_{nc}=(q^x;q^c)_n (q^{x+1};q^c)_n \dots (q^{x+c-1};q^c)_n$$
for $x=-N$ and $x=N+1$, we can identify $S_N$ with the LHS of \eqref{1002}
where the parameters are chosen as follows.
Let $s=(c+1)/2$ and $k=b+s$. Suppose $N=mc+t$ with $0\leq t\leq c-1$.
We consider the limit $b_i, c_i\to \infty$ for $i=s+1,\dots, k$. We put
$b_s=q^{t-N}$ and $c_s=q^{Nc+c}$. For $j=1,2,\dots, s-1$,
among the integers $\{0,1,2,\dots,c-1\}$
there is exactly one $u_j$ and $v_j$, such that
$u_j=j+t\pmod c$ and $v_j=-j+t\pmod c$.
We choose $b_j=q^{u_j-N}$ and $c_j=q^{v_j-N}$ for $j=1,2,\dots, s-1$.
The base $q$ in the identity should be replaced by $q^c$. Therefore,
in the rest of this section
$$(q^a)_m:=(q^a;q^c)_m\, .$$

The RHS of the identity gives us the following expression for $S_N$.
\be
S_N=\sum_{N\geq n_k\geq n_{k-1}\geq \dots\geq n_2\geq 0}
\hat F(n_k,\dots,n_2) \tilde F(n_k,\dots,n_2)
\ee
where
$$ \hat F(n_k,\dots,n_2)\sim
\frac{(q^c)_N (q^{-Nc})_{n_k} (q^{Nc+c})_{n_s} (q^{-mc})_{n_s}
(q^{mc-Nc})_{n_{s+1}-n_s}
}{(q^{-Nc})_{n_{s+1}} (q^{c+mc})_{n_{s+1}}\prod^{k-1}_{i=1}(q^c)_{n_{i+1}-n_i}
} \prod^{s-1}_{j=1}(q^{c+2N-v_j-u_j})_{n_{j+1}-n_j}\, ;
$$
$$\tilde F(n_k,\dots,n_2)\sim \prod^{s-1}_j
\frac{(q^{u_j-N})_{n_j}(q^{v_j-N})_{n_j}
}{(q^{c+N-u_j})_{n_{j+1}} (q^{c+N-v_j})_{n_{j+1}}
}\, .
$$
Here $x\sim y$ means $x/y$ is a unit in $\BZ[q^{\pm 1}]$.
Note that $c-2N-v_j-u_j$, which is equal to $2N-2t\pm c$ or $2N-2t$, is always
a multiple of $c$.

Observe that  $\hat F(n_k,\dots,n_2)\neq 0$
  iff the following inequalities  hold
\be \label{10} n_k\leq N \, , \;\;\;\;\;  n_s\leq \lfloor N/c\rfloor=m \ee
(otherwise $(q^{-Nc})_{n_k}$ or $(q^{-mc})_{n_s}$ is zero);
\be \label{20} n_{s+1}-n_s\leq N-m \ee
(otherwise $(q^{mc-Nc})_{n_{s+1}-n_s}$ is zero).

Let us  assume
that $q$ is a primitive $r$--th root of unity,
then we have in addition
\be \label{25} N\leq r/c\, ,\; \;\;\;\; Nc+cn_s< r \ee
(otherwise $(q^{c})_{N}$ or  $(q^{Nc+c})_{n_s}$ is zero).
Note that if  $\hat F(n_k,\dots,n_2)\neq 0$ then it is also
well--defined.

\begin{lemma} Suppose $q$ is a primitive $r$--th root of unity,
then $z \tilde F(n_k,\dots,n_2)$ is divisible by $\widetilde {(q;q)}_N$.
\end{lemma}
\begin{proof}
It suffice to show that $z$ is divisible by $\widetilde {(q;q)}_N D$,
where $D$ is the denominator of $\tilde F(n_k,\dots,n_2)$.
Since $n_2\leq n_3\leq \dots\leq n_s$, we have
$$D\;|\; (q^{1+N})_{n_s}(q^{2+N})_{n_s}\dots (q^{c+N})_{n_s}=
\widetilde {(q^{1+N};q)}_{cn_s}\, $$
and so $\widetilde {(q;q)}_N D$ divides $\widetilde {(q;q)}_N
\widetilde {(q^{1+N};q)}_{cn_s}=\widetilde {(q;q)}_{N+cn_s}$,
but $N+cn_s\leq (r-1)/2$.
Indeed, $$2N+2cn_s\leq 3N+cn_s \leq Nc+cn_s< r$$
by \eqref{10}, \eqref{25}. Hence,
$$\widetilde {(q;q)}_{N+cn_s}\;\; |\;\; \widetilde {(q;q)}_{(r-1)/2}=z\, .$$
\end{proof}

\begin{lemma} For a primitive $r$--th root of unity $q$,
$\hat F(n_k,\dots,n_2)$ is divisible by $\widehat {(q;q)}_N=
(q^c;q^c)_m$.
\end{lemma}
\begin{proof}
Using for integer $a\geq b>0$ the formula
$$(q^{-ac})_b\sim \frac{(q^c)_a}{(q^c)_{a-b}},$$
we have
$$
\frac{(q^{-Nc})_{n_k}}{(q^{-Nc})_{n_{s+1}}}
\frac{\prod^{s-1}_{j=1}(q^{c+2N-v_j-u_j})_{n_{j+1}-n_j}
}{\prod^{k-1}_{i=1}(q^c)_{n_{i+1}-n_i} }\sim
\frac{(q^{c})_{N-n_{s+1}}}{(q^{c})_{N-n_{k}}}
\frac{\prod^{s-1}_{j=1}(q^{c+2N-v_j-u_j})_{n_{j+1}-n_j}
}{\prod^{k-1}_{i=1}(q^c)_{n_{i+1}-n_i} }.
$$

The latter,  using the fact that $(q^c)_a$ divides
$(q^{c+2N-v_j-u_j})_a$, is divisible by
$\frac{1}{(q^c)_{n_{s+1}-n_s}}$. Thus $\hat F(n_k,\dots,
n_2)/(q^c)_m$ is divisible by
$$\frac{(q^c)_{N-m}}{(q^c)_{n_{s+1}-n_s}}\;
\frac{ (q^c)_{N+n_s}
}{
 (q^c)_{m+n_{s+1}} (q^c)_{N-m-n_{s+1}+n_s}
}\; (q^{-mc})_{n_s}\, .$$
Note that in the first factor  the denominator  divides
the numerator due to \eqref{20}, and in the second factor
because of the binomial integrality.
\end{proof}

\subsection{Diagonalizing the linking matrix}

We say that a closed 3--manifold is of diagonal type if it can be
obtained by integral surgery along an algebraically split link.

\begin{proposition} Suppose $M$ is a closed 3--manifold. There exists
lens spaces $M_1,\dots,M_k$ of the form $L(2^l,a)$ such that the
connected sum of $(M\# M)$ and these lens spaces is of diagonal
type. \label{001}
\end{proposition}

We modify the  proof of a similar result in \cite{Le}.

\subsubsection{Linking pairing}
Recall that {\em linking pairing} on a finite Abelian group $G$ is a
non--singular symmetric bilinear map from $G\times G$ to $\BQ/\BZ$.
Two linking pairing $\nu, \nu'$ on respectively $G,G'$ are
isomorphic if there is an isomorphism between $G$ and $G'$ carrying
$\nu$ to $\nu'$. With the obvious block sum, the set of equivalence
classes of linking pairings is a semigroup.

One type of linking pairing is given by  non--singular square
symmetric matrices with integer entries: any such $n\times n$ matrix
$A$  gives rices to a linking pairing $\phi(A)$ on $G= \BZ^n /A
\BZ^n$ defined by $\phi(A)(v,v') = v^t A^{-1} v' \in \BQ \mod \BZ$,
where $v,v'\in \BZ^n$. If there is a {\em diagonal} matrix $A$ such
that a linking pairing $\nu$ is isomorphic to $\phi(A)$, then we say
that $\nu$ is {\em of diagonal type}.

Another type of pairing is the  pairing $\phi_{b,a}$, with $a,b$
non--zero co--prime integers, defined on the cyclic group $\BZ/b$ by
$\phi_{b,a}(x,y)= axy/b \mod \BZ$. It is clear that $\phi_{b,\pm 1}$
is also of the former type, namely, $\phi_{b,\pm 1}= \phi(\pm b)$,
where $(\pm b)$ is considered as the $1\times 1$ matrix with entry
$\pm b$.

\begin{proposition} Suppose
$\nu$ is a linking pairing on a finite group $G$. There are pairs of
integers $(b_j,a_j),j=1,\dots,s$ with $b_j$ a power of 2 and  $a_j$
either $-1$ or $3$, such that the block sum of $\nu\oplus \nu$ and
all the $\phi_{b_j,a_j}$ is of diagonal type. \label{1101}
\end{proposition}

\begin{proof} The following pairings, in 3 groups, generates the
semigroup of linking pairings, see \cite{Kojima,Wall}:

Group 1: $\phi(\pm p^k)$, where $p$ is a prime, and $k>0$.

Group 2: $\phi_{b,a}$ with $b=p^k$ as in group 1, and $a$ is a
non-quadratic residue modulo $p$ if $p$ is odd, or $a=\pm 3$ if
$p=2$.

Group 3: $E_0^k$  on the group $\BZ/2^k \oplus \BZ/2^k$ with $k \ge
1$ and $E_1^k$ on the group $\BZ/2^k \oplus \BZ/2^k $ with $k \ge
2$.

For explicit formulas of $E_0^k$ and $E_1^k$, see \cite{Kojima}. We
will use only a few relations between these generators, taken from
\cite{Kojima,Wall}.

Any pairing $\phi$ in group 1 is already diagonal by definition,
hence $\phi\oplus \phi$ is also diagonal.

A pairing $\phi=\phi_{b,a}$ in group 2 might not be diagonal, but
its double $\phi \oplus \phi$ is always so: Suppose $b$ is odd, then
one of the relations is $\phi_{b,a}\oplus \phi_{b,a} = \phi(b)
\oplus \phi(b)$, which is diagonal type. Suppose $b$ is even, then
$b=2^k$, $a=\pm 3$, and one of the relations says $\phi_{b,\pm
3}\oplus \phi_{b,\pm 3} = \phi(\mp b) \oplus \phi(\mp b).$

Thus $\nu\oplus \nu$ is the sum of a diagonal linking pairing and
generators of group 3.

 Some of the relations concerning group 3 generators are

 \begin{align*}
 E^k_0 \oplus
\phi_{2^k, -1} &= \phi(2^k)\oplus \phi(-2^k) \oplus \phi(-2^k) \\
E_1^k \oplus \phi_{2^k,3}&= \phi(2^k)\oplus \phi(2^k)\oplus
\phi(2^k).
\end{align*}

Thus by adding to $\nu\oplus \nu$  pairings of the forms $\phi_{2^k,
a}$ with $a=-1$ or $a=3$, we get a new linking pairing which is
diagonal.
\end{proof}

\subsubsection{Proof of Proposition \ref{001}} Every closed
3--manifold $M$ defines a linking pairing, which is the linking
pairing on the torsion of $H_1(M,\BZ)$. Connected sum of 3--manifolds
corresponds to block sum of linking pairings.

First suppose $M$ is a rational homology 3--sphere, i.e.  $M$ is
obtained from $S^3$ by surgery along a framed oriented link $L$,
with non--degenerate linking matrix $A$.
 Then the linking pairing on $H_1(M,\BZ)$ is exactly
$\phi(A)$. Also, the lens space $L(b,a)$ has linking pairing
$\phi_{b,a}$. Proposition \ref{001} follows now from Proposition
\ref{1101} and the well--known fact that
 if the linking pairing on $H_1(M,\BZ)$ is of
diagonal type, then $M$ is of diagonal type, see \cite{Ohtsuki,Le}.

The case when $M$ has higher first Betti number reduces to the case
of rational homology 3--spheres just as in \cite{Le}.

\subsection{Proof of Theorem \ref{main} (general case)}

\begin{lemma} Suppose $(a,r)=1$, and $M=L(a,b)$, then lens space.
Then $\tau_M(\xi) \in \BZ[\xi]$ and moreover, $\tau_M(\xi)$ is
invertible in $\BZ[\xi]$. \label{003}
\end{lemma}
\begin{proof}
 This follows from the explicit formula for the $SO(3)$ invariant of
a lens space given below \eqref{lens}. Note that if  $a^* a=1\pmod r$,
then
$$\frac{1-\xi}{1-\xi^{a^*}}=\frac{1-\xi^{a a^*}}{1-\xi^{a^*}}\, .$$
\end{proof}

\vspace*{0.2cm} \noindent {\bf Proof of Theorem \ref{main} (general
case)} Choose the lens space $M_1, \dots, M_k$ as in Proposition
\ref{001}. Since $N:= M\#M \#M_1\#\dots \#M_k$ is of diagonal type,
its $SO(3)$ invariant is in $\BZ[\xi]$. Note that the orders of the
first homology of $M_1,\dots,M_k$ are  powers of $2$, and hence
co--prime with $r$. Lemma \ref{003} shows that the $SO(3)$ invariant
of $M\#M$ is in $\BZ[\xi]$, and by Lemma \ref{2001},  the
$SO(3)$ invariant of $M$ is in $\BZ[\xi]$ too. \qed

\section{Rational surgery formula}

\subsection{Hopf chain}
Let  $a,b$ be  co--prime integers  with $ b>0$.
It is well known that rational surgery with parameter $a/b$ over a
link component can be achieved by shackling that component
with a framed Hopf chain and then performing integral surgery,
in which the framings $m_{1,2,...,n}$ are related to $a/b$ via a
continued fraction expansion:
$$\frac{a}{b}=-\frac{1}{\displaystyle m_n-\frac{1}{\displaystyle
    m_{n-1}-\dots \frac{1}{\displaystyle m_2-\frac{1}{\displaystyle
        m_1}}}}
$$
Let $D:= (F_{U^{+}}(\xi))^{\sigma^H_+}\,
(F_{U^{-}}(\xi))^{\sigma^H_-}$ where $\sigma^H_\pm$ is the number
of the  (positive/negative) eigenvalues of the linking matrix for
 the Hopf chain.
Let $\left( \frac{d}{r}\right)$ be the Jacobi symbol and $s(b,a)$ the
Dedekind sum. Recall that

$$ s(b,a) := \sum_{i=1}^{|a|-1}\left(\left( \frac{i}{a}
\right)\right) \left(\left( \frac{ib}{a} \right)\right), \quad
\text{where $\left(\left( x \right)\right):= x-\lfloor x \rfloor
-1/2$}.
$$

\begin{lemma} \label{hopf}
For odd $r$ with $(b,r)=1$, we have
$$\frac{[j]}{D}
{\sum_{j_1,...,j_n}}^\xi \prod^{n}_{i=1} q^{m_i\frac{j_i^2-1}{4}} [j_i]
\psdiag{7}{21}{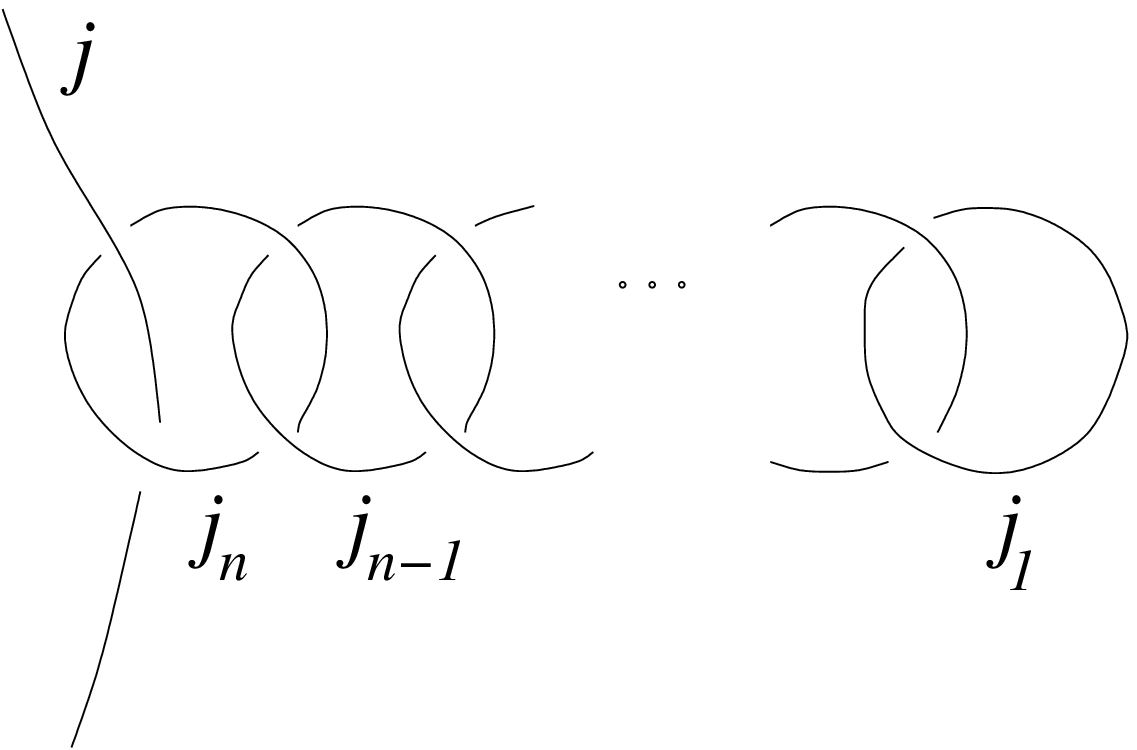}=\left(\frac{b}{r}\right)q^{3s(a,b)}\left[\frac{j}{b}\right]q^{\frac{a(j^2-1)}{4b}
}\psdiag{7}{21}{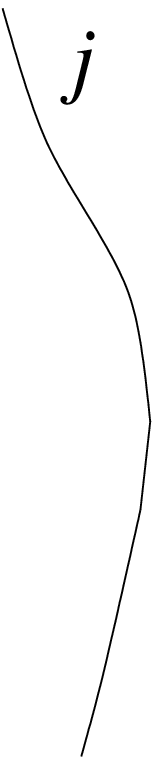}
$$
\end{lemma}

\begin{proof}
The colored Jones polynomial of the $(j_1,j_2)$--colored Hopf link
is $[j_1j_2]$. Thus we have to compute
$${\sum_{j_1,...,j_n}}^\xi q^{\sum_i m_i\frac{j_i^2-1}{4}}
(q^{\frac{j_1}{2}}-q^{-\frac{j_1}{2}})(q^{\frac{j_ij_2}{2}}-q^{-\frac{j_ij_2}{2}})\dots
(q^{\frac{j_n j}{2}}-q^{-\frac{j_nj}{2}})
$$
The result is given by Lemma 4.12 in \cite{LL}, where $A=q^{1/4}$
has the same order as $q$, because $r$ is odd. Moreover, $p$ and $q$
in \cite{LL} are related to our parameters as follows:
$a=-q$ and $b=p$. Computations analogous to Lemmas 4.15--4.21 in \cite{LL}
imply the result.
\end{proof}

If $(r,a)=1$,
the $SO(3)$ invariant of
the lens space $L(a,b)$, which is obtained by surgery along the
unknot with  rational  framing $a/b$, can
be easily computed.
\begin{equation}
 \tau_{L(a,b)}(\xi) = \left(
\frac{a}{r}\right)\,\ev_\xi\left( q^{-3 s(b,a)} \frac{ q^{1/2a} -
q^{-1/2a} }{q^{1/2} - q^{-1/2}}\right).\label{lens}
\end{equation}
Here we used the
 Dedekind reciprocity law (see e.g. \cite{KM}), where  $\sn(d)$
is the sign of $d$,
\be\label{reciproc}
12\left(s(a,b)
  +s(b,a)\right)=\frac{a}{b}+\frac{b}{a}+\frac{1}{ab}-3\sn(ab)\, ,
\ee
 multiplicativity of the Jacobi symbols
$\left(\frac{ab}{r}\right)=\left(
\frac{a}{r}\right)\left(\frac{b}{r}\right)$ and
\begin{equation}
\frac{\gamma_d(\xi)}{\gamma_{\sn(d)}(\xi)} = \left
(\frac{|d|}{r}\right )\,\ev_\xi(q^{(\sn(d)-d)/4}). \label{gamma}
\end{equation}
which holds for any nonzero integer $d$.
Note that
$\tau_{L(a,b)}(\xi)$ is invertible in $\BZ[\xi]$.

\subsection{Laplace transform}
Laplace transform method, developed in \cite{BBL}, allows us
to construct  unified invariant
by  computing the Laplace transform of $\qbinom{n+k}{2k+1} [n]$,
and by   proving its divisibility
by $\frac{\{2k+1\}!}{\{k\}!}$.
Let us show how this strategy works for rational framings.



Suppose that one component of $L$ has rational framing $a/b$.
Then by Lemma \ref{hopf} we have to compute
$$\cL_{a/b;n} \left( \qbinom{n+k}{2k+1} \{k\}!
\left\{\frac{n}{b}\right\}\right) = \frac{\{k\}!}{\{2k+1\}!}
\cL_{a/b;n}\left(\frac{\{n/b\} \, \{n+k\}!}{\{n-k-1\}!}\right)\, .$$

Let $Y_k(q,n,b):=\{n/b\} \, \{n+k\}!/\{n-k-1\}!$.
One can easily see that $Y_k(q,n,b)=Y_k(q,-n,b)$ and
$Y_k(q,n,b)=Y_k(q^{-1},n,b)$.
This implies for  $H_k(q,a/b):=\cL_{-a/b;n}(Y_k(q,n,b))$ that
$$H_k(q,a/b)=H_k(q^{-1},-a/b)\; .$$
Therefore, it is sufficient to compute
$H_k(q,a/b)$ for  $a>0$.

\subsection{Divisibility of the Laplace transform images}

\begin{proposition} \label{110x}

For  $a,b\in \N$, $(a,r)=1$, $(b,r)=1$
and $k\leq (r-3)/2$,  we have
$$ {\sum_{n}}^\xi  q^{\frac{a(1-n^2)}{4b}} \,  \qbinom{n+k}{2k+1} \{k\}!
\left\{\frac{n}{b}\right\}  = 2 \,q^{\frac{(b-1)^2}{4ab}}
\gamma_{-a/b}(\xi) \, \ev_\xi (
F_k(q,a,b))\, .$$
where
$F_k(q,a,b)\in q^{\frac{(3k+2)(k+1)}{4}} \R_{a,{2k+1}}$.
\end{proposition}


A similar formula in the case  $b=1$ was obtained in
\cite{Le}. Proposition \ref{110x} implies that
\be \label{relat} F_k(q,a,1)= \frac{\{k\}!}{\{2k+1\}!}
Y(k,a)\ee
where $Y(k,a)$ was defined in \cite{Le} as follows.
\be
Y(k,a)  :=  \sum_{j=0}^{2k+1}(-1)^j
\qbinom{2k+1}{j}\,
q^{(j-k)^2/a}
\ee

The proof of Proposition \ref{110x} is given in Appendix.
In the rest of the section we define $F_k(q,a,b)$.
We put
$$ C_{k,a,b}=
(-1)^k q^{\frac{(5k+2)(k+1)}{4}} t^{\frac{k(k+1)}{2}(2b-3)}
\frac{(t;t)_{2k+1}}{(q;q)_{2k+1}} $$
Let $w$ be
a primitive root of unity of order $a$.
Let $t$ be the $a$--th primitive root of $q$, i.e. $t^a=q$.
We use the following notation
 $(q^x)_y=(q^x;q)_y$, $(t^x)_y=(t^x;t)_y$ and
$(w^{\pm i} t^y)_x=(w^i t^y;t)_x (w^{-i}t^y;t)_x$.
\noindent
{\bf Case $a$ is odd.}
For { odd} $a$ we define
$c=(a-1)/2$, $l=c+b-1$, $x_i=\sum^{c-1+i}_{j=1} m_j$.

\begin{multline}\label{odd}
\frac{F_k(q,a,b)}{ C_{k,a,b}}:=
\sum_{m_{1}, \dots , m_{l}\ge 0, x_b\leq k}
(-1)^{m_1}
t^{-\frac{m_1(m_1+1)}{2} +\sum^b_{i=1} x_i(x_i-1)}
\frac{(t^{4k+2})^{m_{c-1}+2m_{c-2}+ \dots
+(c-1)m_1}}{(t^{2k})^{m_l+2m_{l-1}+\dots+lm_1}}\\
\frac{(q^{-k})_{k-m_1} (t^{2k+2})_{m_1}(t^{2k+2})_{m_2}\dots (t^{2k+2})_{m_c}}
{(t)_{k-x_b} (t)_{m_2} (t)_{m_3} \dots (t)_{m_l} }
\frac{(w^{\pm 2}t^{-2k-1})_{m_1}
 (w^{\pm 3} t^{-2k-1})_{m_1+m_2}
\dots  (w^{\pm c}t^{-2k-1})_{x_0} }
{(w^{\pm 2} t^{m_1+1})_{m_2}
(w^{\pm 3} t^{m_1+1})_{m_2+m_3}\dots
(w^{\pm c} t^{m_1+1})_{x_1-m_1}}
\end{multline}

\noindent
{\bf Case $a$ is even.}
For even $a$ we define
$c=a/2-1$, $l=c+b$, $x_i=\sum^{c+i}_{j=1} m_j$.

\begin{multline}\label{even}
\frac{F_k(q,a,b)}{ C_{k,a,b}}:=
\sum_{m_{1},  \dots , m_{l}\ge 0, x_b\leq k}
(-1)^{m_1}
t^{-\frac{m_1(m_1+1)}{2} +\sum^b_{i=1} x_i(x_i-1)}
\frac{(t^{4k+2})^{m_{c-1}+2m_{c-2}+ \dots
+(c-1)m_1}
}{(t^{2k})^{m_l+2m_{l-1}+\dots+lm_1}}\\
\frac{(t^{3k+1})^{x_{-1}}(-w^ct^{4k+2})^{m_c}\;\;
(q^{-k})_{k-m_1} (t^{2k+2})_{m_1}\dots (t^{2k+2})_{m_{c-1}}(-w^{-c}
t^{k+1})_{m_c}(-w^{c}t^{2k+2})_{m_{c+1}} }
{(t)_{k-x_b} (t)_{m_2} (t)_{m_3} \dots (t)_{m_l} } \\
\frac{  (w^{\pm 2}t^{-2k-1})_{m_1}
  (w^{\pm 3}t^{-2k-1})_{m_1+m_2}
\dots  (w^{c} t^{-2k-1})_{x_{-1}}(w^{-c} t^{-2k-1})_{x_{0}}
 (-t^{-2k-1})_{x_0} }
{(w^{\pm 2} t^{m_1+1})_{m_2}
(w^{\pm 3} t^{m_1+1})_{m_2+m_3}\dots
(w^c t^{m_1+1})_{x_{0}-m_1} (-t^{-k+x_{-1}})_{m_{c}}
(w^{-c} t^{m_1+1})_{x_1-m_1}(-t^{m_1+1})_{x_1-m_1}}
\end{multline}

\noindent
{\bf Example.}
\be\label{IHS}
F_k(q,1,b) := q^{-\frac{(3k-2)(k+1)}{4}} q^{kb(k+1)}
\sum_{m_{1}, m_{2}, \dots , m_{b}\ge
0,\;\sum m_i=k}
 q^{\sum^{b-1}_{i=1}(x_i^2-(2k+1)x_i)}
\frac{(q)_{k}}
{\prod_{i=1}^{b}
(q)_{m_i}} \ee
where $x_p=\sum^p_{i=1} m_i$.

Note that (\ref{IHS}) coincides up to units with
the formula for the
coefficient $c'_{k,b}$ in the decomposition of $\omega^{b}$
computed in \cite[(46)]{Mas}. (The same coefficient
(up to units) appears  also in the cyclotomic expansion
of the Jones polynomial of twist knots (\ref{twist})). This is because
 surgery on the $(-1/b)$--framed component can be achieved by replacing
this component by $b$ 1--framed copies. Indeed,
changing the variables in (\ref{IHS}) as follows:
$s_1=k-m_1$, $s_2=k-m_1-m_2,...,s_{b-1}=k-m_1-\dots-m_{b-1}$, we get
$$F_k(q,1,b)=q^{\frac{(k+2)(k+1)}{4}}
\sum_{k\geq s_1\geq s_2\geq \dots\geq s_{b-1}\geq 0}
q^{s^2_1+s^2_2+\dots +s^2_{b_1}+s_1+\dots+s_{b-1}}
\frac{(q)_k}{\prod^{b-1}_{i=1}(q)_{s_i-s_{i+1}}}\, .
$$

\section{Universal invariant}
In this section  we   assume that $(r,a)=1$, where
$r$ is the order of the root of unity $\xi$ and $a=|H_1(M,\BZ)|$.

Let  $M=L(a,b)$ be a lens space with $a>0$. Then the universal invariant
$I_M$ was defined in \cite{Le} as follows.
$$ I_{M} :=    q^{3s(1,a)-3 s(b,a)}
\frac{1-q^{-1/a}}{1-q^{-1}}.$$
Note that $3(s(1,a)-s(b,a))\in \Z$ and $I_M$ is invertible in $\Gd$.

For an arbitrary rational homology sphere $M$
with $a=|H_1(M,\BZ)|$, it was shown in \cite{Le}
that there are lens spaces $M_1,\dots,M_l$ such that
$M'= (\#_{i=1}^l M_i) \# M$ can be obtained by surgery on an algebraically
split link and  $I_{M_j}$
are  invertible in $\Gd$. Then we can define
$$ I_{M} = I_{M'} ( \prod_{i=1}^l I_{M_i})^{-1}.$$

It remains to define $I_M$
when $M$ is given by surgery along an algebraically split link $L$.
Assume $L$ has
 $m$ components with nonzero rational framings  $\frac{a_1}{b_1}, \dots
\frac{a_m}{b_m}$. Then we have    $|H_1(M,\Z)|=a$ for $a=\prod_i a_i$.
Let $L_0$ be $L$ with all framings switched to zero.

\begin{theorem}
 For $M$ as above,  the unified  invariant
is given by the following formula.
\begin{equation}\label{univ1}
I_M= q^{(a-1)/4}\sum_{k_i=0}^{\infty}
J_{L_0}(P'_{k_1},\dots,P'_{k_m}) \prod_{i=1}^m \sn(a_i)
q^{\frac{1}{2a_i} -3s(b_i,a_i)}
F_{k_i} (q^{-\sn(a_i)},|a_i|,b_i)
\end{equation}
Moreover,

$$\left(\frac{a}{r}\right) \tau_M(\xi)=\ev_\xi (q^{(1-a)/4} I_M)$$
\end{theorem}

\begin{proof}
First observe that,
if $b_i=1$ for all $i$, our formula coincides with (21) in \cite{Le}.
It follows from  (\ref{relat}) and
$$q^{\frac{3\sn(a_i)-a_i}{4}} q^{3s(1,a_i)}=q^{\frac{1}{2a_i}}\; .$$
 Here we used that
$3s(1,a)=\frac{1}{2a}
+\frac{a-3{\rm sn}(a)}{4}$ by the reciprocity law (\ref{reciproc}).

Let us collect the coefficients in the definition of $\tau_M$.
From Lemmas \ref{un1}, \ref{hopf}, Proposition
 \ref{110x} and (\ref{gamma}) we have
$$q^{3s(a_i,b_i)-\frac{(b_i-1)^2}{4a_ib_i}+ \frac{3{\sn}(a_i)}{4}-\frac{a_i}{4b_i}}
=q^{-3s(b_i,a_i) } q^{\frac{1}{2a_i}}$$

The Corollary 0.3 (d) in  \cite{Le} allows to drop
the conditions $(b_i,r)=1$, because
$I_M$ is determined by its values at any infinite sequence of roots of
prime power order  from $U_a$.
\end{proof}

\subsection{Proof of  Theorem \ref{ring} }
The statement holds trivially if $M=L(a,b)$. Indeed, we have $m=1$,
$f_{k}=0$ for $k>0$, and $f_0=q^{3s(1,a)-3 s(b,a)} (1-t)/(1-q)\in \R_1$.

The general case follows from  \eqref{univ1}
and Proposition \ref{110x}. Note that multiplication of
$I_M$ by the inverse of $I_{L(a,b)}$ multiply all  $f_{k_i}$
by an element of $\BZ[t^{\pm 1}]$.
 Moreover,
 $I_M$ does not contain fractional powers of $q^{1/a}$
(compare  Proof of Lemma 4.2 in \cite{Le}).
\qed

\section{Applications}

In this section we compute  the universal invariant
$I_M$ for Seifert fibered spaces and for $a/b$ surgeries on twist knots.

\subsection{Seifert fibered spaces with a spherical base}
Let $M=L(b;a_1/b_1,\dots,a_n/b_n)$ be the Seifert fibered space
with base space $S^2$ and with $n$ exceptional fibers with
orbit invariants $(a_i,b_i)$ ($a_i>0, 0\leq b_i\leq a_i$, $(a_i, b_i)=1$),
and with bundle invariant $b\in \Z$.

It is well--known that $M$ is a rational homology sphere
if $e:=b+\sum b_i/a_i\neq 0$ and $|H_1(M,\Z)|=|e|\prod_i a_i$.
Furthermore, $M$ can be obtained by surgery on the following
(rationally framed) link.

$$\psdiag{7}{21}{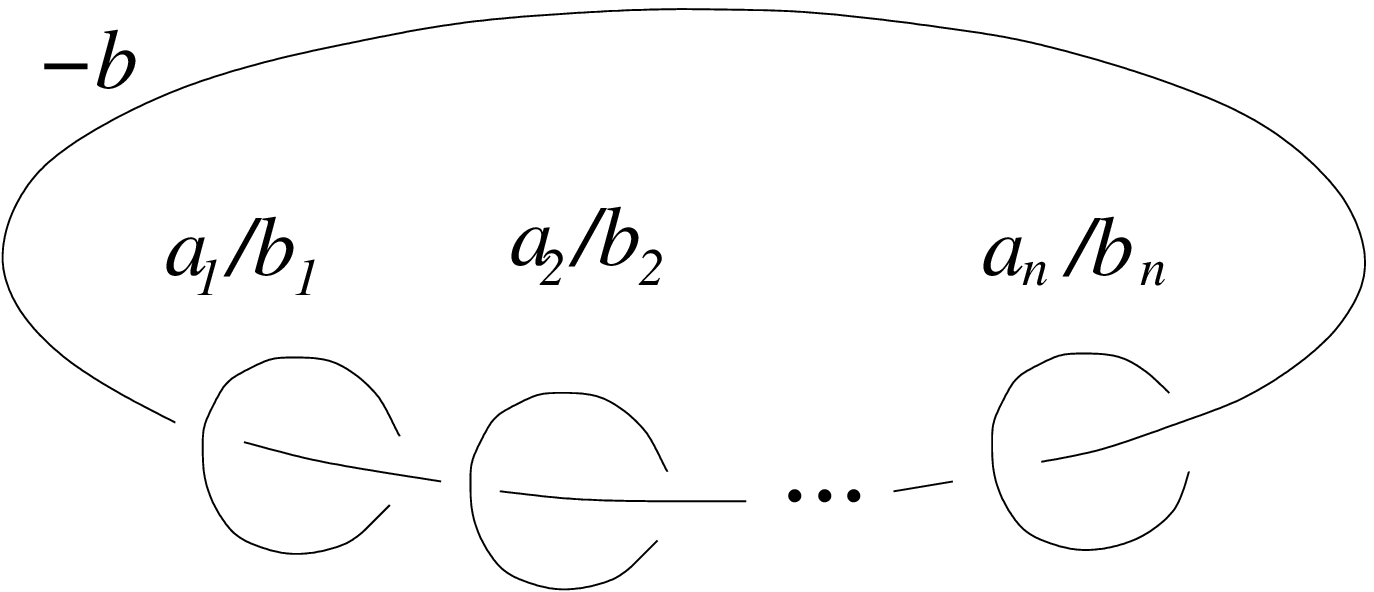}$$
\begin{theorem} Let  $M=L(b;a_1/b_1,\dots,a_n/b_n)$ as above.
Assume $e\neq 0$,  and $|H_1(M,\Z)|=d$.

$$I_M= q^{\frac{d-1}{4}}\;
\frac{q^{(e-3\sn(e))/4}q^{-3\sum_i s(b_i,a_i)} }{\{1\}}
\cL_{-e;j}
\left(\frac{\prod^n_{i=1}\left\{ \frac{j}{a_i}
\right\}
}{\left\{ {j}\right\}^{n-2}}  \right )
$$

\end{theorem}

\begin{proof}
The linking matrix of the surgery link has $n$ positive eigenvalues
and the sign of the last eigenvalue is equal to  $-\sn(e)=-\sn(b)$.
Let us color the rationally framed components of the surgery link by
$j_i$, $i=1,\dots, n$ and the $-b$ framed component by $j$.

The main ingredient of the proof is the following computation.
Using  Lemmas \ref{un1},  \ref{hopf} we have
$$
\left(\frac{b_i}{r}\right)
\frac{q^{3s(a_i,b_i)}}{F_{U^+(\xi)}}
\sum^\xi q^{\frac{a_i(j_i^2-1)}{4b_i}}
\left[\frac{j_i}{b_i} \right]
[jj_i]=\left(\frac{a_i}{r}\right)
q^{-3s(b_i,a_i)}q^{-\frac{b_i(j^2-1)}{4a_i}}\left[\frac{j}{a_i} \right]
$$

Applying finally  the Laplace transform  $\cL_{-e;j}$
  and collecting the factors we get the result.
\end{proof}

\subsection{Proof of Theorem \ref{SFS}}
Note that
$M=L(b;a_1/b_1,\dots,a_n/b_n)$ is an integral homology sphere if $e^{-1}
=\pm\prod_i a_i$. $M$ is uniquely determined by
the pairwise co--prime integers $a_i$. (Knowing $a_i$'s and $e$,
one can compute  $b_i$'s and $b$ using the Chinese remainder theorem).

Suppose for simplicity that $e>0$.
Rewriting
$$\frac{1}{\left\{ {j}\right\}^{n-2}} =(-1)^{n-2}q^{(n-2)/2}
\sum_{k=0} c_k q^k$$
with $c_k\in \Z$, we see that the image of the Laplace transform
is the sum of the following terms:
$$(-1)^{n-2}c_k\, q^{\frac{\prod_i a_i}{4}(\pm 1/a_1\pm 1/a_2\dots\pm1/a_n
+2k+n-2)^2} $$
The leading term   in $I_M$  for $k\to \infty$ behaves
asymptotically like
$$q^{k^2\prod_i a_i  + k(n-2)\prod_i a_i +
k\sum_i a_1\dots \hat a_i \dots a_n }
$$
where $\hat a_i$ means delete $a_i$.
This allows to determine the $a_i$'s. In the case $e<0$,
we have the same asymptotic after replacing $q$ by $q^{-1}$.
$\hfill\Box$
\vspace*{.1cm}

Theorem \ref{SFS} follows also from \cite{BL}, where the
computation were done for the LMO invariant combined with the $sl_2$ weight
system, i.e. for  the Ohtsuki series. By the result of Habiro,
$I_M$ is as powerful as the  Ohtsuki series.

\subsection{Dehn surgeries on twist knots}

Let $K_p$ be the twist knot with $p$ twists.
Masbaum \cite{Mas} calculated the $P'_n$ colored Jones polynomial of this knot.
For $p>0$ we have
\be\label{twist}
J_{K_p}(P'_n)= q^{n(n+3)/2}\sum_{i_1,i_2,\dots,i_p\geq 0, \sum_j i_j=p}
q^{\sum_i (s^2_i+s_i)}\frac{(q)_{n}}{\prod^p_{j=1}(q)_{i_j}}\ee
where  $s_k=\sum^k_{j=1} i_j$.
The formula for the negative $p$ can be obtained from the given one
by sending $p\to -p$, $q\to q^{-1}$, forgetting the factor $q^{n(n+3)/2}$
and multiplying the result by $(-1)^n$.

\begin{corollary}
Let $M_{a/b}$ is obtained by $(a/b)$ surgery in $S^3$
on the twist knot $K_p$.
Then
\begin{equation}
I_{M_{a/b}}:= q^{(a-1)/4} \sn(a)
q^{-3s(b,a)+\frac{1}{2a}}
\sum_{n=0}^{\infty}
J_{K_p}(P'_n)
F_n (q^{-\sn(a)},|a|,b)
\end{equation}
\end{corollary}

\subsection{Proof of Theorem \ref{detector}}
Assume  $K$ and $K'$ are  0--framed.
We expand the function
  $Q_{K} (N):= J_K(N)[N]$ around $q=1$ into power series.
 Suppose $q=e^h$, then we have
$$Q_{K}(N)|_{q=e^h}=\sum_{2j\leq n+2}
c_{j,n}(K)  N^j h^n\, .$$ It is known that $c_{j,n}$ is zero if $j$
is odd. Applying  Laplace transform,
 we have to replace $N^{2j}$ by $\frac{(-2)^j(2j-1)!! }{a^j h^j}$ (see \cite{Le2}).
Therefore, the following expression coincides (up to some standard
factor) with the Ohtsuki series
$$\sum_{2j\leq n+2} c_{2j,n}(K) (-2)^j(2j-1)!! \; a^{-j}h^{n-j}\, .$$
From $I_{M(K,a)}=I_{M(K',a)}$ we derive
$$\sum_{2j\leq n+2}
\left(c_{2j,n}(K)-c_{2j,n}(K')\right) (-2)^j(2j-1)!!\; a^{-j}
h^{n-j}=0\, .$$  Because the last system of equations should hold
for infinitely many $a\in\BZ$, we have $c_{2j,n}(K)=c_{2j, n}(K')$
and $J_K(N)=J_{K'}(N)$ for any $ N\in \N$.
\qed

\section*{Appendix}

The main technical ingredient we  use in the proof of Proposition
\ref{110x}  is the Andrew's generalization of the  Watson's identity
(\cite[Theorem 4, p.199]{And}):

\begin{multline}\label{Watson}
_{2p+4}\phi_{2p+3}\left[\begin{array}{l}
\alpha,t\sqrt{\alpha},-t\sqrt{\alpha},b_1,c_1,...,b_p,c_p,t^{-N}\\
\sqrt{\alpha},-\sqrt{\alpha},\alpha t/b_1,\alpha t/c_1,...,\alpha
t/b_p,\alpha t/c_p, \alpha t^{N+1}
\end{array}\; ; t, \frac{\alpha^p t^{p+N}}{b_1c_1...b_pc_p}\right]=
\frac{(\alpha t)_N(\alpha t/b_pc_p)_N}{(\alpha t/b_p)_N(\alpha t/c_p)_N}$$
\\
 \sum_{m_1,...,m_{p-1}\geq 0}
\frac{(b_p)_{\sum_i m_i}(c_p)_{\sum_i m_i}
(t^{-N})_{\sum_i m_i}}{(b_pc_pt^{-N}/\alpha)_{\sum_i m_i} }
\prod^{p-1}_{i=1} \frac{t^{ m_i}(\alpha t)^{(p-i-1)m_i} (\alpha t/b_ic_i)_{m_i}
(b_i)_{\sum_{j<i} m_j} (c_i)_{\sum_{j<i}m_j}
}{(t)_{m_i}
(\alpha t/b_i)_{\sum_{j\leq i}m_j}(\alpha t/c_i)_{\sum_{j\leq i}m_j}
(b_i c_i)^{\sum_{j< i}m_j}
}
\end{multline}
where
\begin{equation}\label{Wa}
_{r}\phi_{s}\left[\begin{array}{l} a_1, a_2,\dots, a_r\\
b_1,\dots,b_s\end{array}; t,z\right]=\sum^{\infty}_{n=0}
\frac{(a_1)_n(a_2)_n\dots(a_r)_n}{(t)_n(b_1)_n\dots
  (b_s)_n} \left[ (-1)^nt^{\binom{n}{2}}\right]^{1+s-r} z^n
\end{equation}
are the basic q--hypergeometric series and $(a)_n:=(a;t)_n$.
\subsection*{Proof of Proposition \ref{110x}}

We have to compute
 $\cL_{-a/b;n}(\{n/b\} \, \{n+k\}!/\{n-k-1\}!)$.
Note
$$\{n/b\} \, \{n+k\}!/\{n-k-1\}!)=
q^{-n/2-nk-n/(2b)}(1-q^{n/b})(q^{n-k})_{2k+1}\, .$$
Using the q--binomial theorem and \eqref{four} (with $c=1$) we get
$$
 q^{\frac{(2bk+b+1)^2}{4ba}}\sum^{\infty}_{j=0}
\frac{(q^{-2k-1})_j}{(q)_j} q^{\frac{b}{a}j^2+(1-2\frac{b}{a})kj+
(1-\frac{(b+1)}{a})j}(1-q^{\frac{2j-2k-1}{a}})$$
We put $t:=q^{1/a}$ and choose a primitive $a$--th root of unity $w$.
Then using
$$(q^x;q)_j= (t^x;t)_j(wt^x;t)_j(w^2 t^x;t)_j\dots
(w^{a-1}t^x;t)_j$$
we can rewrite the previous sum as follows.
\begin{equation} \label{LHS}
\sum^{2k+1}_{j=0}\frac{(t^{-2k-1})_j(wt^{-2k-1})_j\dots
(w^{a-1}t^{-2k-1})_j}{(t)_j(wt)_j\dots  (w^{a-1}t)_j}
t^{bj^2+(a-2b)kj+(a-b-1)j}(1-t^{2j-2k-1})
\end{equation}
The main point is that (\ref{LHS}) is equal to $(1-t^{-2k-1})$ times
the LHS of the generalized Watson identity (\ref{Watson})  with the
specialization of parameters described below. We  consider  the limit
 $\alpha\to t^{-2k-1}$.
We set $p={\rm max}\{b,a+b-2\}$,  $b_i,c_i\to \infty$ for $i=a-1,...,p-1$;
$b_p\to t^{-k}$, $c_p\to \infty$ and
$N\to\infty$.

\noindent
{\bf Case $a$ is odd.}
We put $c=(a-1)/2$;
$b_i=w^i t^{-2k-1}$, $c_i=w^{-i}t^{-2k-1}$ for $i=1,...,c$;
$b_i, c_i\to t^{-k}$ for $i=c+1,...,a-2$.

\noindent
{\bf Case $a$ is even.}
We put $c=a/2-1$. Let $p=a+b-2$,
$b_i=w^i t^{-2k-1}$, $c_i=w^{-i}t^{-2k-1}$ for $i=1,...,c-1$;
$b_{c}=w^c t^{-2k-1}$, $c_{c}=-t^{-k}$;
$b_{c+1}=-t^{-2k-1}$, $c_{c+1}=w^{-c}t^{-2k-1}$;
$b_i, c_i\to t^{-k}$ for $i=c+2,...,a-2$.

To simplify
 (\ref{Watson})
we use the following limits.
\begin{align*}
\lim_{c\to \infty} \frac{(c)_n}{c^n} & =(-1)^n t^{n(n-1)/2}&
\lim_{c\to \infty} \left(\frac{t}{c}\right)_n & =1\\
\lim_{c_1,c_2\to
\infty}\frac{(c_1)_n(c_2)_n}{(t^{-N}c_1c_2)_n}&=(-1)^n t^{n(n-1)/2}
t^{Nn} & \lim_{\alpha\to t^{-2k-1}}\frac{(\alpha t)_{\infty}}{(\sqrt{\alpha t})_{\infty}}
& =2 (t^{-2k})_k
\\
\end{align*}


Finally, the formulas below allow us to separate the factor
$\frac{(t)_{2k+1}}{(q)_{2k+1}}$.
\begin{align*}
\frac{\{2k+1\}!}{\{k\}!}&=q^{-\frac{(3k+2)(k+1)}{4}}(-1)^{k+1}(q^{k+1})_{k+1}\\
(q)_j&=(-1)^{(k-j)}q^{\frac{(j-k)(k+j+1)}{2}} \frac{(q)_k}{(q^{-k})_{k-j}}\\
(t^{-k})_j& =(-1)^j t^{-kj+\frac{j(j-1)}{2}}\frac{(t)_k}{(t)_{k-j}}
\end{align*}

The next lemma implies the result.
$\hfill\Box$

\begin{lemma}
$\frac{F_k(q,a,b)}{C_{k,a,b}}\in \Z[t^{\pm 1}]\, .$
\end{lemma}

\begin{proof}
First note that $F_k(q,a,b)$ does not depend on $w$, because in the LHS
of the identity $w$ does not occur.

Suppose $a$ is odd. Let $z:=x_1-m_1=m_2+m_3+\dots+m_c$.
Let us complete $(w^c t^{m_1+1})_{z}$
to $(q^{m_1+1})_{z}$ by multiplying the numerator and the
denominator of  (\ref{odd}) with
$$ (t^{m_1+1})_z (w^\pm t^{m_1+1})_z (w^{\pm 2} t^{m_1+m_2})_{z-m_2}
 \dots (w^{\pm
(c-1)} t^{m_1+z-m_c})_{m_c}\; .$$
Now up to units the denominator of (\ref{odd}) is equal to
${(q^{m_1+1})_{x_1-m_1} (t)_{k-x_b} (t)_{m_2} (t)_{m_3} \dots (t)_{m_l} }$
which  divides the numerator.
The even case is similar.
\end{proof}

\ifx\undefined\bysame
    \newcommand{\bysame}{\leavevmode\hbox
to3em{\hrulefill}\,}
\fi

\end{document}